\documentclass[a4paper]{article}
\usepackage{amsmath}
\usepackage{amsfonts}
\usepackage{amssymb}
\usepackage{amsthm}
\usepackage{authblk}

\newtheorem{theorem}{Theorem}
\newtheorem{lemma}{Lemma}
\newtheorem{proposition}{Proposition}

\newtheorem{definition}{Definition}

\newtheorem{corollary}{Corollary}

\begin{document}
\title{Binary strings of finite VC dimension}
\author{Hunter Johnson}
\affil{\small Department of Mathematics and Computer Science, John Jay College, CUNY, New York, New York 10019 \\ \textit{E-mail address}, \texttt{hujohnson@jjay.cuny.edu}}

\date{Jan 16, 2020}
\maketitle
\begin{abstract}

Any binary string can be associated with a unary predicate $P$ on $\mathbb{N}$.  In this paper we investigate subsets named by a predicate $P$ such that the relation $P(x+y)$ has finite VC dimension.  This provides a measure of complexity for binary strings with different properties than the standard string complexity function (based on diversity of substrings).  We prove that strings of bounded VC dimension are meagre in the topology of the reals, provide simple rules for bounding the VC dimension of a string, and show that the bi-infinite strings of VC dimension $d$ are a non-sofic shift space. Additionally we characterize the irreducible strings of low VC dimension (0,1 and 2), and provide connections to mathematical logic. 

\end{abstract}

\section{Introduction}

By ``string'' we mean a sequence of symbols from some alphabet $\mathcal{A}$.  There are many natural examples of strings.  This document is a string in the Roman alphabet.  We conceptualize DNA and computer code as strings either in the alphabet of nucleotides, or ones and zeros.  Some strings are more complicated than others. For example the infinite decimal expansion of the number 0 is simple. In genetics, genes often occur in complex regions of the genetic code, whereas other regions of DNA may be monotonous or repetitious.

Formally, if $s$ is a sequence of symbols from $\mathcal{A}$, let $p_s(n)$ denote, for a positive integer $n$, the number of distinct substrings of $s$ of length $n$.  This assumes that the length of $s$ is at least $n$.  By basic combinatorics, 

\begin{equation}
	1 \leq p_s(n) \leq |\mathcal{A}|^n
	\label{eqn:basiccomplexity}
\end{equation}

The function $p_s(\cdot)$ as described above is a commonly used measure of complexity, with applications in dynamical systems, automata theory, biology, and other domains \cite{berthe2010combinatorics}. It features as a principal object in some well known theorems, such as the Morse-Hedlund theorem (see Section \ref{complexVC}.)

There is also a tradition of measuring the complexity of set families.  For example characterizing the families of events for which frequency converges uniformly to probability (across the family) has been of interest in statistics since the 1970's \cite{dudley2014uniform}.  This viewpoint has found applications in the theoretical foundations of machine learning by way of the Vapnik Chervonenkis dimension.  Surprisingly a somewhat parallel investigation in mathematical logic examining uniformly definable families has turned up interesting connections to the ``wildness'' of the semantics of a formal theory \cite{simon2015guide}.

The formal description of Vapnik Chervonenkis dimension is as follows. Consider a set $U$ and a collection of subsets $\mathcal{C} \subseteq \{c : c \subseteq U\}$.  For a given $B \subseteq U$, we define $\mathcal{C} \cap B = \{c\cap B: c \in \mathcal{C}\}$.  The following is a complexity function somewhat analogous to the string complexity function:

\begin{equation}
	m_{\mathcal{C}}(n) = \sup \{|\mathcal{C} \cap B| : B\subseteq U \text{ and } |B|=n\} 
	\label{eqn:vccomplexity}
\end{equation}

It was independently discovered by Sauer, Vapnik and Chervonenkis, and Shelah and Perles that for any $\mathcal{C}$, either $m_{\mathcal{C}}(n)$ always equal to its maximum possible value of $2^n$, or else the function is bounded by a polynomial in $n$.  For a given $\mathcal{C}$, the largest $n$ for which $m_{\mathcal{C}}(n) = 2^n$ is known as its Vapnik-Chervonenkis dimension (if no such $n$ exists, the VC dimension is infinite).  It is standard to refer to $B$ as \textit{shattered} if $\mathcal{C}\cap B$ equals the power set of $B$.  Observe that the VC dimension of $\mathcal{C}$ is the size of the largest shattered subset of $U$ \cite{sauer1972density,shelah1972combinatorial,vapnik2015uniform}.  If the VC dimension of $\mathcal{C}$ is $d$, then $m_{\mathcal{C}}(n) \leq n^d$ for all $n$. 

Valiant discovered that, in the context of a learning problem, hypothesis spaces with finite VC dimension coincide with the hypothesis spaces that are distribution-free learnable (in the Probably-Approximately-Correct model of learning) \cite{valiant1984theory}.  In model theory, Shelah, Laskowski and many others have investigated first order theories in which all partitioned formulas have finite VC dimension \cite{laskowski1992vapnik}. The idea of the VC dimension of a binary string is in some ways implicit in model theory, as explained in Section \ref{Smt}.  However this is not necessarily obvious, and the emphasis here is not on the model theoretic properties of structures, only in strings of finite VC dimension. 

On the combinatorics side, many achievements have been made in characterizing the complexity of arrangements of geometrical objects in a way that essentially relates to VC dimension, for example Radon, Cover, Basu and many others \cite{matousek2013lectures}. Work has also been done on concept classes relating to the sets of positivity for neural networks and other sophisticated learning machines (for example Sontag, Macintyre and many others \cite{sontag1998vc,karpinski1997polynomial}). However nothing seems to have been written about VC classes corresponding to binary strings.  

In model theory there has been work done that relates directly to this paper.  In particular many authors have made deep contributions to understanding when the structure $(\mathbb{N},+,P)$ or $(\mathbb{Z},+,P)$ is stable or NIP when $P$ is a unary predicate (please see Section \ref{Smt} for attributions and definitions). A model theorist would understand this paper as investigating the relation $P(x+y)$ in and of itself, independent of the stability theoretic properties of the structure in which it lies.  Therefore this is not a model theory paper, and the questions we answer are not discussed in the model theoretic literature.  Namely determining VC dimension for particular strings in heuristic ways, dimension preserving alterations of strings, the topological and dynamical nature of strings of finite VC dimension when considered in aggregate, and what we call here ``prime strings" which are combinatorially extreme vis a vis their length and VC dimension.

\subsection{The structure of this paper}


In Section \ref{Sintro} we give the central definitions, including the VC dimension of finite strings. In Section \ref{Sexamples} we give a number of examples of strings of both finite and infinite VC dimension. We investigate a string that shares a self-similarity quality with the Cantor set, and show that it lies on the extreme end of VC complexity. We also investigate strings that are more like Sidon sets (on the other end of the complexity spectrum).  It arises that $p_s(n)$ gives surprisingly little information on the VC dimension of a string (and \textit{vice-versa}). We conjecture that Sturmian strings (which are of minimal non-trivial complexity with respect to $p_s(n)$) can have infinite VC dimension. 

 In Section \ref{Svccomp}, we first study real numbers with binary representations of finite VC dimension.  Theorem \ref{Tcantorlike} establishes that real numbers of uniformly bounded VC dimension are distributed in a Cantor set like fashion; we also show that there are uncountably many such numbers. We investigate the set of reals of finite VC dimension and establish its basic topological properties in Theorem \ref{Ttop}.  We conjecture that the reals of finite VC dimension constitute an uncountable subfield of $\mathbb{R}$, analogous to the constructable numbers or computable numbers.
 
 We then move on to a more detailed examination of binary strings of finite VC dimension.  Theorem \ref{Talt} establishes that the VC dimension of a string can be bounded in terms of the number of alternations between 0 and 1, echoing many related results in model theory. We then investigate the strings of finite VC dimension as a dynamical system and show that the bi-infinite strings of uniformly bounded VC dimension are a non-sofic shift space (Theorem \ref{Tsofic}). 

Then we turn to investigate so-called prime strings, which are strings of a given VC dimension that are atomic in a certain sense. We completely characterize the prime strings of low dimension (0,1, and 2) and prove that prime strings of dimension 3 are fundamentally more complex. Along the way we show that prime strings do not constitute a regular language in the hierarchy of formal languages. 

The last section ties this paper to work in  model theory. 

\section{From strings to set families}\label{Sintro}

In this section we explain how to encode a binary string as a set family, and consequently equip binary strings with a VC dimension.  This involves investigating three different but similar notions of dimensionality.  We then give an application of these definitions to a particular string, and analyze the relationship between VC complexity (of a string) and the complexity as usually defined.

Let $S$ be a binary string (a string in the alphabet $\mathcal{A} = \{0,1\}$), and $\mathcal{S}(S)$ the set of its finite substrings.  If $S$ is infinite, then by default we understand $S$ to be infinite only on the right (we will discuss bi-infinite strings in Section \ref{Sbi}). Recall that substrings, unlike subsequences, are contiguous.  Let $s \in \mathcal{S}(S)$.  We let $s_i$ denote the digit in the $i$th position of $s$ (indexing is zero based).  We will associate $s$ with a subset of natural numbers $n(s) = \{i: s_i = 1\}$. Using this technique we can associate $S$ with a subset family $\mathfrak{S} = \{n(s):s \in \mathcal{S}(S)\}$.  From here we can define the VC dimension of $\mathfrak{S}$ in the usual way. By abuse of notation, the VC dimension of a binary string S refers to the VC dimension of $\mathfrak{S}$.

\textbf{Example}:

If $S = 011$, $\mathcal{S}(S) = \{<>,0,1,01,11,011\}$ and $\mathfrak{S} = \{\{\},\{0\},\{1\},\{0,1\},\{1,2\}\}$.  The set $\{0,1\}$ is shattered and VC($\mathfrak{S}$) = 2.

This way of assigning VC dimension to a binary string is perhaps the most natural, though there are other notions.  Note that any concept class (ie set system) can be regarded as consisting of binary strings, and a concept class derived from a binary string can be thought of as a special kind of concept class that can be assembled to form a single object (ie the string).

We now consider a ``sliding window'' method of assigning a VC dimension to a string.  Given a  string $S$, let $\mathcal{S}_w(S)$ denote the substrings of $S$ of length $w$.  We can imagine sliding a width $w$ window along the string and recording which substrings are observed. $\mathcal{S}_w(S) = \{s \in \mathcal{S}(S): |s| = w\}$.  This leads to its own set family $\mathfrak{S}_w = \{n(s) : s \in \mathcal{S}_w(S)\}$.

\textbf{Example}:

If $S = 011$, $\mathcal{S}_2(S) = \{01,11\}$ and $\mathfrak{S}_2 = \{\{1\},\{0,1\}\}$.  

The sliding window dimension of $S$ is defined to be $SWdim(S) = \max \{VCdim(\mathfrak{S}_w): w \in \mathbb{N} \}$.  If no maximum exists then we say $SWdim(S) = \infty$.

  Often we only care whether a certain measure of complexity is finite or infinite.  From this point of view it is unlikely to greatly matter exactly which notion of string complexity is selected.  We show below that the notions of SWdim and VCdim for binary strings differ in value by at most one.

\begin{proposition} For any binary string $S$,  $SWdim(S) \leq VCdim(S)$.
\end{proposition}

\begin{proof}
 For a given $w$, $\mathcal{S}_w(S) \subseteq \mathcal{S}(S)$.  Therefore for each $w$, $VCdim(\mathfrak{S}_w) \leq VCdim(\mathfrak{S})$. Thus $\max \{VCdim(\mathfrak{S}_w) : w \in \mathbb{N}\} \leq VCdim(\mathfrak{S})$. 
\end{proof} 

\begin{proposition} For any binary string, if $VCdim(S) \geq d+1$ then $SWdim(S) \geq d$.  

\end{proposition}
\begin{proof} Suppose $\mathfrak{S}$ shatters the set of indices $\{i_0,...,i_d\}$.  Let $\mathfrak{S}^{i_d} = \{n(s) : s \in \mathcal{S}(S), s_{i_d} = 1\}$. 

Then $\mathfrak{S}^{i_d}$ shatters $\{i_0,...,i_{d-1}\}$.  Moreover every $s$ in $\mathcal{S}(S)$ such that $s_{i_d} = 1$ must have $|s| \geq {i_d}+1$.

Then $\mathfrak{S}^{i_d}_{i_d+1}$ also shatters $\{i_0,...,i_{d-1}\}$.

Therefore $\mathfrak{S}_{i_d+1}$ shatters $\{i_0,...,i_{d-1}\}$, and this proves $SWdim(S) \geq d$. \end{proof}

The lower bound in Corollary \ref{C1} is seen to be tight by the examples provided above.  The upper bound can be established using $S=01100$.

\begin{corollary}\label{C1}
The above two propositions together show that $$VCdim(S)-1 \leq SWdim(S) \leq VCdim(S).$$

In particular, $SWdim(S)=\infty \iff VCdim(S)=\infty$.
\end{corollary}

To ease the transition to further generalizations, we introduce one additional notion of string complexity.  

Given a string $S$, a $d$-mask on $S$ is a set of subsequences of $S$ of the form $\langle S_{i_0+t},S_{i_1+t},\ldots,S_{i_{d-1}+t}\rangle$, where $t$ varies over the index set (in this case the natural numbers). The indexes must be distinct -- for convenience we can assume $i_0 < i_1 < \cdots < i_{d-1}$. We ignore sequences resulting from values of $t$ that give nonsensical indexes (too big or too small).  

A $d$-mask is said to be \textit{full} if it contains all binary sequences of length $d$. 

\textbf{Example}:

$S = 101001$.

Then one 2-mask might be $\{\langle S_{0+t},S_{2+t}\rangle : t = 0,1,2,3\}$.  Explicitly this is $\{\langle 11\rangle,\langle 00\rangle,\langle 10\rangle,\langle 01 \rangle\}$.  This is a full 2-mask.

The mask dimension $Mdim(S) =\max \{d : S\text{ has a full $d$-mask}\}$.  If there is no maximum then we say $Mdim(S) = \infty$. 

\begin{proposition}\label{SWeqM} For all binary strings $S$, $Mdim(S) = SWdim(S)$.
\end{proposition}

\begin{proof} 

Suppose there is a full $d$-mask of the form $\{\langle S_{i_0+t},S_{i_1+t},\ldots,S_{i_{d-1}+t} \rangle : t \in \mathbb{N}\}$.  Without loss we may assume that $i_0 = 0$.  Let $w=i_{d-1}+1$ (assuming $i_0 < i_1 < \cdots < i_{d-1}$). Then for every $t$ there is some $s$ in $\mathcal{S}(S)_w$ such that $n(s) \cap \{i_0,i_1,\ldots,i_{d-1}\} = \{i_j : S_{i_j+t}=1\}$. Because the $d$-mask is full, $\mathfrak{S}_w$ shatters $\{i_0,i_1,\ldots,i_{d-1}\}$.

Conversely suppose $\mathfrak{S}_w$ shatters $\{i_0,i_1,\ldots,i_{d-1}\}$. For any $d$ length binary sequence there is some $s\in \mathcal{S}(S)_w$ such that $\langle s_{i_0},s_{i_1},\ldots,s_{i_{d-1}} \rangle$ realizes the sequence. If $t$ is the starting index of $s$ in $S$ then $\langle S_{i_0+t},S_{i_1+t},\ldots,S_{i_{d-1}+t}\rangle = \langle s_{i_0},s_{i_1},\ldots,s_{i_{d-1}} \rangle$. Thus $\langle S_{i_0+t},S_{i_1+t},\ldots,S_{i_{d-1}+t}\rangle$ where $t$ varies over the index set gives a full $d$-mask.

\end{proof}

The following is easy, but interesting; if a string is infinitely complex iff it has an infinitely complex tail. 

\begin{lemma}\label{Ltail}  A binary string S has finite VC dimension iff it has a suffix of finite VC dimension.
\end{lemma}
\begin{proof}  From left to right is obvious. For right to left, let $S = xy$ where $y$ is a suffix of finite VC dimension $d$.  If $S$ has infinite VC dimension then it must have arbitrarily large mask dimension.  Therefore there must be a full $|x|+d+1$ mask.  The last $d+1$ coordinates of this mask are a full $d+1$ mask.  However all parts of $S$ affecting the substrings contributing to the full $d+1$ mask must occur in $y$.  This is a contradiction.
\end{proof}

\section{Examples and comparisons with $p_s(n)$}\label{Sexamples}

Having labored through the above definitions and propositions, we now explore in some depth a particular example.  

\subsection{The Cantor string}
We inductively define a string which is somewhat like the Cantor set.  We call this construction the Cantor string. 
Let $0^k$ denote a string of $k$ zeros for $k \in \mathbb{N}$.  

Let $S^{(0)} = 1$.  Given $S^{(i)}$ for $i \in \mathbb{N}$, define $S^{(i+1)} = S^{(i)}0^{|S^{(i)}|}S^{(i)}$.  That is, $S^{(i+1)}$ is the concatenation of $S^{(i)}$, $|S^{(i)}|$ many zeros, and $S^{(i)}$ again.
The first few examples of $S^{(i)}$ are given below.

$S^{(0)} = 1$ 

$S^{(1)} = 101$ 

$S^{(2)} = 101000101 $

$S^{(3)} = 101000101000000000101000101 $

$S^{(4)} = 101000101000000000101000101000000000000000000000000000101000101000000000101000101 $

Note that $S^{(i)}$ is a proper initial segment of $S^{(i+1)}$ for all $i$.  The \textit{Cantor string}, $S^{(\omega)}$ is the infinite binary string that has all $S^{(i)}, i \in \mathbb{N}$, as a proper initial segment.

\begin{lemma}\label{cantorindexes}  The indexes in $S^{(\omega)}$ where a 1 occurs are precisely the indexes of the form $2\cdot\sum_{k\in K}3^k$ for some finite $K \subseteq \mathbb{N}$ (possibly empty).
\end{lemma}
\begin{proof}  This is an easy inductive argument on $i$ for the $S^{(i)}$.
\end{proof}

Just as the Cantor set shows that sparsity and cardinality are independent, the Cantor string shows that a sparse string can have high complexity. 

\begin{proposition}  The Cantor string has infinite VC dimension.
\end{proposition}
\begin{proof}  It suffices to show that the string has infinite mask dimension.  Let $d \in \mathbb{N}$ be given.  We will show that $S^{(\omega)}$ admits a full $d$-mask.  The indexes we will use are of the form $i_l := 2\cdot 3^l$ for $l=1,2,\cdots,d$. Let $s$ be a finite proper initial substring of $S^{(\omega)}$ which is long enough that all the indexes referenced below exist.

Let $[d]$ denote the set $\{1,2,\ldots,d\}$.  Let $A \subseteq [d]$ be given.  We will construct a value $t$ such that for all $l \in [d]$, $s_{2(3^l)+t} = 1 \iff l \notin A$.  

Define $t = 2(\sum_{a\in A} 3^a)$.  Consider $2(3^l)+t = 2(3^l+\sum_{a \in A}3^a)$.  Observe that $3^l+\sum_{a \in A}3^a$ is a sum of cubic powers iff $l \notin A$.  

Thus, by Lemma \ref{cantorindexes}, $s_{2(3^l)+t} = 1 \iff l \notin A$.  This means that we have constructed a full $d$-mask.  Since $d$ was arbitrary, the dimension is infinite. 

\end{proof}

\subsubsection{The complexity function and the Cantor string}\label{complexVC}

We say that a string is \textit{aperiodic} if it is not periodic, meaning it does not factor as the infinite product of some word. A string is \textit{not eventually periodic} if no suffix of the string is periodic. 

Consider the following well-known result \cite{morse1938symbolic}.

\begin{lemma}\label{LMHT}[Morse-Hedlund theorem]  An aperiodic sequence has a strictly increasing complexity function.
\end{lemma}

The authors of the above theorem characterized the simplest strings that are not eventually periodic as being Sturmian.  A string $S$ is \textit{Sturmian} if $p_s(n) = n+1$ for all $n$.

The Cantor string is clearly aperiodic, since it contains arbitrarily long substrings of the form $ 10^{3k}1$.  For this same reason, the Cantor string is not even eventually periodic.  It is also unbalanced, in the sense that the Hamming weight of substrings of length $n$ take on multiple (more than two) values for large $n$. Sturmian strings are characterized by being balanced and not eventually periodic, and so the Cantor string is not Sturmian.  However the complexity function is still linearly bounded.

\begin{theorem} When $x$ is the Cantor string, $p_x(n) = 2n-1$ for all $n>1$.
\end{theorem}

\begin{proof}  
  We use the strings $S^{(k)}$ as in the definition of the Cantor string.  First observe that $p_x(1) = 2$. We claim that for $n \geq 2$, $p_x(n) = 2n-1$. Let $k$ be maximal so that $n>3^{k-1}$. Because of the recursive structure of the Cantor string, in order to determine $p_x(n)$ it suffices to consider $p_w(n)$ where $w=0^{n}S^{(k)}0^{n-1}$.  Recall that $S^{(k)} = S^{(k-1)}0^{3^{k-1}}S^{(k-1)}$.  We will consider the distinct words taken on by a window of length $n$ as it progresses across $w$.  More precisely, we imagine that we have a sequence of binary $n$-tuples denoted by $b$, where $b(t)=\langle w_{0+t},w_{1+t},\ldots,w_{n-1+t}\rangle$ as $t$ ranges over $\{0,1,2,\ldots,n+3^k-1\}$.

First note $b(0)=0^n$.  This and all subsequent words are distinct until $t=n+3^{k-1}$.  At this stage, graphically, the leftmost index of $b(t)$ is positioned at the start of the middle 0's in $S^{(k)}$.  We then encounter words previously seen for the next $2\cdot 3^{k-1}-n+3^{k-1}$ values of $t$ (incrementing sequentially).  Graphically, at this stage, the rightmost index of $b(t)$ is moving into the long field of zeros on the right after $S^{(k)}$.  Then the next $n-(3^{k-1}+1)$ values of $t$ again yield words that are previously unseen.  

The total number of steps performed is $n+3^k-1$, which is one step for each value of $t$.  The total number of distinct words encountered is $n+3^{k-1}+n-(3^{k-1}+1) = 2n-1$.  

\end{proof}

The apparent contradiction between infinite VC dimension and tame string complexity will be addressed in Section \ref{Svccomp}.

\subsection{The Thue-Morse sequence}

By the \textit{parity function} we mean the function $f:\mathbb{N} \rightarrow \{0,1\}$ such that $f(n)=0$ iff the binary representation of $n$ has an even number of 1's. The Thue-Morse sequence \cite{allouche1999ubiquitous} is the sequence for which the $n$th symbol is $f(n)$. Here we show that the Thue-Morse sequence has infinite mask dimension (and hence infinite VC dimension).  We must show that the sequence admits a full $d$-mask for all $d$. Let $d \in \mathbb{N}$ be given.  Construct a binary matrix $M$ of dimensions $d \times 2^{d+1}$ in the following way. We assume some canonical mapping $g : \mathcal{P}([d]) \rightarrow \{1,2,...,2^d\}$, where $\mathcal{P}([d])$ is the powerset of $\{1,2,3,\ldots,d\}$.

We conceive of $M$ as being composed of $d$ rows and $2^d$ two-digit columns.  Then for $A \in \mathcal{P}([d])$, $2g(A)$ will index the first (ie. leftmost) digit of the column corresponding to $A$, and $2g(A)+1$ will index the second digit. For each $i \in [d]$, the row $i$ of $M$ will have a 1 precisely in the columns $\{2g(A)+1 : i \in A\}$.

We now conceptualize the rows of $M$ as integers $a_1,a_2,\ldots,a_{d}$ where the row $i$ of $M$ specifies the binary digits of $a_i$.

Let $A \subseteq [d]$ be given.  Let $b_A = 2^{2g(A)+1}$, the integer whose binary expression is 1 precisely in index $2g(A)+1$.

Then $f(b_A + a_i) = 0 \iff i \in A$. This gives a full $d$-mask on the Thue-Morse Sequence.  Therefore the Thue-Morse sequence has infinite mask (and hence VC) dimension.

\subsection{Powers of two and Golumb rulers}\label{Spowersof2}\label{Ssidon}

The facts given in this section relate to recent results in the model theory literature, and nothing in this section is essentially new.  Please see Section \ref{Smodelcon} for a discussion. Sidon sets were explicitly used in the context of VC dimension by \cite{aschenbrenner2016vapnik} and the sequel.

Let $S$ be the infinite binary string with a 1 precisely in indexes that are powers of two.

$$S = 011010001000000010000000000000001\ldots$$

Observe that the index distance between any two 1's is of the form $2^k-2^l = \sum_{k>j\geq l} 2^j$.

Thus if $2^k-2^l = 2^m-2^n$ then (viewing the numbers in binary) it is clear that $k=m$ and $l=n$.

This implies that there can be no full 3-mask.  If there were such a mask $\langle S_{i_0+t},S_{i_1+t},S_{i_2+t}\rangle$ then there must be:

\begin{enumerate}
\item  Some $t_1$ such that $\langle S_{i_0+t_1},S_{i_1+t_1},S_{i_2+t_1}\rangle  = \langle 1,1,1\rangle$
and
\item  Some $t_2$ such that $\langle S_{i_0+t_2},S_{i_1+t_2},S_{i_2+t_2}\rangle = \langle 1,1,0\rangle $
\end{enumerate}

Clearly $t_1 \neq t_2$. But by the above discussion there is at most a unique $t_1$ such that $S_{i_0+t_1}=S_{i_1+t_1}=1$.  Therefore there is no full 3-mask.

On the other hand there is a full 2-mask (namely any two adjacent indexes). Therefore the mask dimension is 2 (which is also the VC dimension).

The $S$ above, when interpreted as the binary expansion of a real number, is a Fredholm constant.  This shows that the binary decimal representation of transcendental numbers can have finite VC dimension.

The proof of the VC dimension of powers of two essentially used only the property that the difference between any pair of 1's is unique.  There is a general term for strings with this property, namely Golumb rulers (or Sidon sets) \cite{o2004complete}.

Let $A \subset \mathbb{N}$. Let $\Delta_A(d) = |\{(a,b)\in A^2: b-a=d, b>a \}|$. We say that $A$ is a Sidon set if $\max_{d \in \mathbb{N}} \Delta_A(d) = 1$.  We will also say that $A$ is a nearly Sidon set if $\max_{d \in \mathbb{N}} \Delta_A(d)$ exists (ie is finite). Any binary string $S$ can be understood as a subset of $\mathbb{N}$, namely $n(S) = \{i : S_i = 1\}$.  We say $S$ is (nearly) Sidon if $n(S)$ is a (nearly) Sidon set.

$S$ is said to be eventually (nearly) Sidon spaced if $S$ is infinite and has a (nearly) Sidon spaced suffix.  

\begin{proposition}
Any of the following conditions imply that a binary string $S$ is nearly Sidon spaced.
\begin{enumerate}
\item $S$ is Sidon spaced.
\item $S$ is eventually Sidon spaced.
\item $S$ is eventually nearly Sidon spaced.
\end{enumerate}
\end{proposition}
\begin{proof}
This fact is straightforward.
\end{proof}
The following are examples of Sidon spaced sequences:

\begin{enumerate}

\item The binary string $S$ with 1 precisely in positions $\{q^n : n \in \mathbb{N}\}$ for a positive integer $q>1$.
\item The binary string $S$ with 1 precisely in positions $\{n! : n \in \omega\}$.

\item Any increasing sequence $\langle u_i : i\in \mathbb{N}\rangle$ such that $u_{i+1} \geq 2u_i$ for all $i$.
\end{enumerate}

\begin{proposition}\label{P:sidon}
Any Sidon spaced binary string $S$ based on a sufficiently large Sidon set has VC dimension 2. 
\end{proposition}
\begin{proof}
This is essentially the same as the proof for powers of two.
\end{proof}

 In fact any nearly Sidon spaced sequence has finite VC dimension. The reason is that for $S$ to have a full $d$ mask implies $\Delta_{n(S)}(k)\geq 2^{d-2}$ for at least one $k \in \mathbb{N}$.  Thus if $S$ has infinite VC dimension then $S$ cannot be nearly Sidon spaced.

\subsection{The number of Sidon sets in $\mathbb{N}$}

A natural question is to ask the cardinality of all strings of a certain VC dimension.  In this subsection we prepare to answer this question, which is formally resolved in Theorem \ref{T:uncountable}.  We will basically show that all strings of any VC dimension can be encoded as strings of VC dimension 2. 

Consider that for binary strings in general, flipping 1's to 0's can increase VC dimension.  For example the simple string 1* can be made to have any desired VC dimension by introducing 0's at certain indexes.  However for Sidon sets this is not the case, because subsets of Sidon sets are also Sidon sets.  Removing elements introduces no new differences or sums, and so the property of being (nearly) Sidon is preserved. 

This property provides an easy way to see that there are uncountably many Sidon sets. Suppose $S$ is a binary string that has 1 precisely at indexes in the set $n(S)=\{i:S_i=1\}$. Let $2^S$ denote the binary string for which $n(2^S)=\{2^i:S_i=1\}$.  In other words $2^S$ has a 1 in index $i$ if and only if $i=2^j$ and $S_j=1$.

\begin{proposition}
For any binary string $S$, $VCdim(2^S)\leq 2$. 
\end{proposition}
\begin{proof}
The string $2^S$ is Sidon spaced.
\end{proof}

The above proposition is interesting, because $2^S$ is in some sense just as complex as $S$.  However this complexity is not something that we can capture with the relation $P(x+y)$ (in the language of Section \ref{Smodelcon}.)  On the other hand the VC dimension of the relation $P(2^{x+y})$ would be capable of detecting ``exponential level'' complexity in $2^S$. There is no reason to stop at first powers, and this discussion could go on to towers of exponentiation or other fast growing functions. From a certain viewpoint, complexity depends on the expressive power of the observer.

\begin{proposition}
For any binary strings $S$ and $T$, $2^S=2^T \iff S=T$. 
\end{proposition}
\begin{proof}
This is obvious.
\end{proof}
  
\begin{corollary}\label{Cuncountable}
There are uncountably many Sidon spaced sequences.
\end{corollary}
\begin{proof}
Let $2^{\mathbb{N}}$ represent the space of all possible binary strings.  This set is uncountable. Then by the previous proposition $\{2^S: S\in 2^\mathbb{N}\}$ is also uncountable.  All of these are Sidon spaced.
\end{proof}

\section{VC dimension and substring diversity}\label{Svccomp}

The Cantor string example shows that for a string to have finite VC dimension it is not sufficient that the complexity function be polynomially (or even linearly) bounded.  However the complexity function $p_S(n)$ can give information about VC dimension in extreme cases.  For example if the complexity function is superpolynomial (on the one hand) or constantly bounded (on the other) then the VC dimension of the associated sequence is determined as either infinite or finite (respectively).  When we establish the finitude of VC dimension for strings in this section, the results are implicit in some model theoretic work, for example \cite{POINT1}. 

\begin{lemma}\label{Lspol}  Suppose that the complexity function $p_S(n)$ for a string S is superpolynomial.  Then S has infinite VC dimension.
\end{lemma}
\begin{proof}
Suppose by way of contradiction that $VCdim(S)=d$ for an integer $d$.  Then by Sauer's Lemma, for any integer $w$, $|\mathcal{S}_w(S)| = |\mathfrak{S}_w|  \leq w^d$. However $p_S(w) = |\mathcal{S}_w(S)|$.
\end{proof}

\begin{corollary} Suppose that $S$ is a binary string of finite VC dimension and let $\mathfrak{L}=\mathcal{S}(S)$ denote the formal language consisting of substrings of $S$. Then $\mathfrak{L}$ is a sparse language.
\end{corollary}
\begin{proof}
This is just a rephrasing of the above Lemma. 
\end{proof}

\begin{lemma}\label{Lperiod} A periodic string $S$ with period $m \in \mathbb{N}$ has mask dimension at most $m$.\end{lemma}
\begin{proof}
Let $\{\langle S_{i_0+t},S_{i_1+t},\ldots,S_{i_{d-1}+t} \rangle : t \in \mathbb{N}\}$ be a full $d$-mask for some integer $d$. Suppose that $i_p \equiv i_q \mod m$ for some $p,q \in \{0,\ldots,d-1\}$.  Then $i_p + t \equiv i_q + t\mod m$, and thus $S_{i_p+t}=S_{i_q+t}$ for all $t$.  Because the $d$-mask is full, we must have $i_p$ distinct modulo $m$ for all $p \in \{0,\ldots,d-1\}$.  Thus $d \leq m$. 
\end{proof}

\begin{corollary} We can draw the following conclusions from Lemma \ref{Lperiod}.

\begin{enumerate}
\item If $S$ has period $m$ then $VCdim(S) \leq m+1$.
\item If $S$ is eventually periodic, then $VCdim(S)$ is finite. 
\item If $p_s(n)$ is bounded by a constant then $VCdim(S)$ is finite.
\item The binary representation of any rational number has finite VC dimension.
\end{enumerate}
\end{corollary}
\begin{proof}
The statement (1) follows from Corollary \ref{C1} and Proposition \ref{SWeqM}.  The statement (2) follows from Lemma  \ref{Ltail}.  The statement (3) follows from Lemma \ref{LMHT}.
\end{proof}

There is a very natural question that we have not answered:

\textbf{Question}:  Can a Sturmian string have infinite VC dimension?

One would think that the answer to this question is ``yes.''  The Cantor string suggests that complexity of a string can grow very slowly provided that the radius of shattered sets of size $d$ grows exponentially in $d$.  A priori we could do this even with minimal non-trivial complexity. But we haven't managed to determine the VC dimension of any Sturmian string to date.

\subsection{The VC dimension of real numbers}

We first state some obvious facts for convenience of reference.

\begin{lemma}\label{Lmonotone}
Let $S$ be a binary string and $S'$ a substring of $S$.  Then $VCdim(S) \geq VCdim(S')$.
\end{lemma}

\begin{lemma}\label{Lfinite}  Let $S$ be an infinite binary string of VC dimension $d$ for $d \in \mathbb{N}$. Then $S$ has a proper initial substring of the same VC dimension.
\end{lemma}

\begin{lemma}\label{Laddzero} If $S$ is a finite binary string and $0^*$ is the zero string, then $VCdim(S0^*)=VCdim(S)$.
\end{lemma}

We now make a few remarks on strings that arise as binary representations of reals.  Rational numbers of the form $\frac{p}{2^k}$ for integers $p,k$ may have more than one base 2 representation. For example the real number 1 has the binary representations $10^*$ and $01^*$.  In this case we will choose the representation that ends in $0^*$.  For the purposes of computing VC dimension, we will ignore any decimal point or negative sign.

A number is said to be normal in base $b$ if, for every positive integer $n$, all possible $n$-digit substrings have density $b^{-n}$ . 

\begin{corollary}  For any real number normal to base 2, the base 2 representation has infinite VC dimension.
\end{corollary}

Borel proved \cite{ian1993borel} that the set of real numbers normal to every base has full Lebesgue measure.  Therefore almost all real numbers have infinite VC dimension with respect to their base 2 representations.  This implies that the real numbers of finite VC dimension are measure zero.  

It is conjectured that for any algebraic irrational, the complexity function is $p_x(n) = O(b^n)$, where $b$ is the base of the representation \cite{adamczewski2007complexity}.  If this conjecture is true then by Lemma \ref{Lspol} the binary decimal representation of any algebraic irrational has infinite VC dimension.

For a string $S$ and nonnegative integers $i<j$, let $S^{([i:j])}$ denote the substring of $S$ defined by $\langle S_i,S_{i+1},\ldots,S_{j-1}\rangle$.  We also write $S^{([i:\infty])}$ for the suffix of $S$ beginning at index $i$ (inclusive).  If we write $S^{([i:-j])}$ this means the same as $S^{([i:n])}$ where $n=len(S)-j$.  If $S$ and $T$ are binary strings, define $L_i^{T,S} =S^{([0:i])}T^{([i:\infty])}$.  The \textit{straight line} from $T$ to $S$ is the sequence of strings $\mathfrak{L}(T,S) = \{L_i^{T,S}:i \in \mathbb{N}\}$.  Define $P_i^{T,S} = S^{([0:i])}T$.  The \textit{push} from $T$ to $S$ is $\mathfrak{P}(T,S) = \{P_i^{T,S}: i \in \mathbb{N}\}$.  If we refer to a limit of $P_i$ or $L_i$ we are referring to the limit of the corresponding real numbers (via binary representation).  For a binary sequence $S$ we let $\mathfrak{r}(S)$ denote the corresponding real number in [0,1). By $0^*$ we mean the binary string with $0$ at all indexes. 

\begin{proposition}\label{Pshoveline}
 Let $S,T$ be binary strings with VC dimension $d_S,d_T \in \mathbb{N}\cup\{\infty\}$ (respectively). The following hold.  
\begin{enumerate}
\item For every $i \in \mathbb{N}$, $VCdim(P_i^{T,S})\geq d_T$.
\item For every $i \in \mathbb{N}$, $VCdim(L_i^{0^*,S})\leq d_S$.
\item $\lim_{i \to \infty} P_i^{T,S} = \mathfrak{r}(S)$
\item $\lim_{i \to \infty} L_i^{T,S} = \mathfrak{r}(S)$
\end{enumerate} 
\end{proposition}
\begin{proof}
The statement (1) follows from Lemma \ref{Lmonotone}. The statement (2) follows from lemmas \ref{Lmonotone} and \ref{Laddzero}.  Statements (3) and (4) are obvious.
\end{proof}

We now consider what happens in the limit to set of real numbers with various VC complexity assumptions. Interestingly, VC dimension can drop arbitrarily in a limit, but not go up.  The following propositions make this precise. 

\begin{proposition}
For any $d \in \mathbb{N}\cup \{\infty \}$ there is a sequence $\langle r_i : i \in \mathbb{N}\rangle$ of real numbers such that $VCdim(r_i)\geq d$ for all $i$ and $\lim_{i \to \infty} r_i = 0$. 
\end{proposition}
\begin{proof}
Let $r_0$ be a real number of VC dimension $d$ and $S$ its binary representation. Let $r_i$ be the real number whose binary representation is $0^iS$ for $i>0$.  Then $\lim_{i\to\infty} r_i = 0$ and $VCdim(r_i)\geq d$ for all $i$. 
\end{proof}

\begin{proposition}\label{Plimits}
For any $d \in \mathbb{N}$ and any sequence $\langle r_i : i \in \mathbb{N}\rangle$ of real numbers such that $\forall i\, VCdim(r_i) \leq d$, if $r^* = \lim_{i \to \infty} r_i$ then $VCdim(r^*)\leq d$. 
\end{proposition}
\begin{proof}
Let $S$ denote the binary representation of $r^*$.  Suppose, by way of contradiction, that $VCdim(S) > d$. Let $S'$ be a finite initial substring of $S$ such that $VCdim(S')>d$.  Let $n = len(S')$.  There is some $i$ such that $|r_i-r^*| < \frac{1}{2^n}$.  Then $S'$ is also an initial substring of the binary representation of $r_i$.  But then $VCdim(r_i)>d$ $\rightarrow \leftarrow$.
\end{proof}

For $d\in \mathbb{N}$, Let $\mathfrak{V}_{\leq d} = \{r\in \mathbb{R}: VCdim(r) \leq d\}$.  Let $\mathfrak{V}_{< \infty}$ denote the reals of finite VC dimension.

\begin{proposition}\label{Pperf}
Let $r \in \mathbb{R}$.  Then the following hold.
\begin{enumerate}
\item There is a sequence of real numbers in $\mathfrak{V}_{< \infty}$ that converges to $r$.  Additionally, if $d=VCdim(r)$ is an integer, then there is a sequence of real numbers in $\mathfrak{V}_{\leq d}$ that converges to $r$.

\item There is a sequence of real numbers in $\mathbb{R}\setminus \mathfrak{V}_{< \infty}$ that converges to $r$.
\end{enumerate}
\end{proposition}
\begin{proof}
Let $T$ be such that $r=\mathfrak{r}(T)$.  Then $\mathfrak{L}(0^*,T)$ witnesses the statement (1).

Let $S$ be any binary string of infinite VC dimension.  Then $\mathfrak{P}(S,T)$ witnesses the statement (2).
\end{proof}

The following theorem shows that $\mathfrak{V}_{\leq d}$ is topologically very similar to the Cantor set.

\begin{theorem}\label{Tcantorlike}
  For all $d\in\mathbb{N}$,  $\mathfrak{V}_{\leq d}$ is closed, totally separated, perfect, and nowhere dense in $\mathbb{R}$. It also has Lebesgue measure 0. 
\end{theorem}
\begin{proof}
By Proposition \ref{Plimits}, every convergent sequence in $\mathfrak{V}_{\leq d}$ has a limit in $\mathfrak{V}_{\leq d}$.  Therefore it is closed.  By  Proposition \ref{Pperf}, every $r\in\mathfrak{V}_{\leq d}$ is a limit point for a sequence in $\mathfrak{V}_{\leq d}$.  Therefore $\mathfrak{V}_{\leq d}$ is perfect.  Given any $r \in \mathfrak{V}_{\leq d}$ let $S$ be its binary representation and take $n > d$ (possibly $n=\infty$).  For any proper initial substring $S'$ of $S$ there is a binary word $w$ such that $S'w$ has VC dimension $n$. This shows that for any open $U$ with $r \in U$, there are real numbers of VC dimension $n>d$.  If $x,y \in \mathfrak{V}_{\leq d}$ are distinct, there is some $z$ strictly between them with $z \notin \mathfrak{V}_{\leq d}$. Thus $\mathfrak{V}_{\leq d}$ is totally separated as witnessed by the open sets $(-\infty,z)\cup(z,\infty)$.

To see that $\mathfrak{V}_{\leq d}$ is nowhere dense, let $U \subseteq \mathbb{R}$ be an open set. Let $r \in U$ be arbitrary. By \label{Pperf} (2) there is a sequence of reals of VC dimension $> d$ converging to $r$.  Therefore there is some  $t\in U\setminus \mathfrak{V}_{\leq d}$. Since $\mathfrak{V}_{\leq d}$ is closed, there is an open set $V$ such that $t\in V$ and $\mathfrak{V}_{\leq d} \cap V=\emptyset$. Then $V\cap U$ shows that $\mathfrak{V}_{\leq d}$ is not dense in $U$.

The theorem of Borel on normal numbers, already mentioned, shows that $\mathfrak{V}_{\leq d}$ has Lebesgue measure 0.

\end{proof}

\begin{corollary}
For any $d\in \mathbb{N}$, set  $\mathfrak{V}_{\leq d} \cap [0,1]$ is a Stone space.
\end{corollary}
\begin{proof}
The set is closed and bounded, therefore compact.  It is obviously Hausdorff.  It is totally separated by Theorem \ref{Tcantorlike}.
\end{proof}

\begin{lemma}\label{Lvc1} If $S$ is a binary string and $VCdim(S)=1$ then $S$ is of one of the following forms. Each of these may be followed by arbitrarily many zeros, which do not affect the VC dimension.  
\begin{enumerate}
\item $0^a1$, $a\geq 1$
\item $1^a$, $a \geq 1$ (possibly $a$ is infinite)
\item $0^a(10^b)^c1$, $0 \leq a \leq b$, $b\geq 1$, $c\geq 1$ (possibly $c$ is infinite).
\end{enumerate}
\end{lemma}
\begin{proof}
The VC dimension is one in each case by easy inspection. In form (3) note that if $a > b$ then the VC dimension is 2. Observe that any nonzero binary string has a string of form (1) or (2) as a prefix.  We will show that any attempt to extend one of these forms while preserving VC dimension results in another instance of one of the three given forms. 

First consider extending a string of form (1). Suppose $S=0^a1$ for $a\geq 1$. Then $VCdim(S1)=2$. The extension $S0^b$ for a nonnegative integer $b$ is trivial.  The extension $S0^b1$ results in a string of type (3) if $b\geq a$, and $VCdim(S0^b1)=2$ otherwise.  We deal with the type (3) possibility below.

Consider extending a string of form (2). Suppose $S=1^a$ for $a \geq 1$. The only nontrivial extension for this form is $S0^b1$ for $b \geq 1$, and $VCdim(S0^b1)=2$.

Now consider extending a string of form (3). Suppose $S=0^a(10^b)^c1$, where $0 \leq a < b$, and $c\geq 1$.  Then $VCdim(S1)=2$. If we let $T=S0^d1=0^a(10^b)^c10^d1$, then the tail of $T$ will be of the form $10^b10^d1.$  But to avoid increasing the VC dimension, we must have $b=d$, yielding another string of form (3).
\end{proof}


\begin{theorem}\label{T:uncountable}
The set $\mathfrak{V}_{\leq d}$ is uncountable iff $d\geq 2$.
\end{theorem}
\begin{proof}
The right to left direction follows easily from Corollary \ref{Cuncountable}. To show countability for $\mathfrak{V}_{\leq 1}$, we first note that there is a unique string of VC dimension 0, namely $0^*$.  Lemma \ref{Lvc1} completes the argument by showing that the strings of VC dimension 1 are of finitely many forms, each with countably many instances.
%
\end{proof}

\begin{theorem}\label{Ttop}
 The following all hold in the standard topology on $\mathbb{R}$.

\begin{enumerate}
\item $\mathfrak{V}_{< \infty}$ is not closed.
\item $\mathfrak{V}_{< \infty}$ is not open.
\item $\mathfrak{V}_{< \infty}$ is dense and codense in $\mathbb{R}$.
\item $\mathfrak{V}_{< \infty}$ is uncountable and co-uncountable in $\mathbb{R}$.
\item $\mathfrak{V}_{< \infty}$ has Lebesgue measure 0.
\item $\mathfrak{V}_{< \infty}$ is meagre.

\end{enumerate}
\end{theorem}

\begin{proof}
 Every real (of any VC dimension) is, on the one hand, the limit of reals of finite VC dimension, and, on the other, the limit of reals of infinite VC dimension.   Therefore the closure of $\mathfrak{V}_{< \infty}$ and the closure of its complement both yield $\mathbb{R}$.  Therefore $\mathfrak{V}_{< \infty}$ is both dense and codense.   Neither $\mathfrak{V}_{< \infty}$ nor its complement are closed under limits, so neither is closed (or open) in $\mathbb{R}$.  Because $\mathfrak{V}_{< \infty} = \bigcup_{d\in\mathbb{N}} \mathfrak{V}_{\leq d}$ and each $\mathfrak{V}_{\leq d}$, is nowhere dense, $\mathfrak{V}_{< \infty}$ is meagre in $\mathbb{R}$.  Because $\mathfrak{V}_{\leq 2}$ is uncountable, $\mathfrak{V}_{< \infty}$ is uncountable. The fact that it is null follows from the theorem of Borel previously mentioned.

\end{proof}



\textbf{Question}:  Is $\mathfrak{V}_{< \infty}$ a subfield of $\mathbb{R}$, intermediate between $\mathbb{Q}$ and $\mathbb{R}$? This subfield of ``simple" reals would be analogous to the constructable or computable numbers. But unlike those examples, it would be uncountable. It is not hard to show that $(\mathfrak{V}_{< \infty},\oplus)$ is a subgroup of $(\mathbb{R},\oplus)$, where $\oplus$ denotes XOR on binary representations. However the ``carry bits" seem to leave just enough uncertainty for the sum of two elements in $\mathfrak{V}_{< \infty}$ to violate closure.  
%
%
%
%
%

\subsection{The structure of strings of finite VC dimension}\label{Sbi}

There is a structural consequence to the topological discussion from the previous section. We present this as Proposition \ref{Pstruct}.  Then we prove some  facts about strings of finite VC dimension, some of which relate to symbolic dynamics. In particular we work toward Theorem \ref{Talt} which shows that VC dimension is bounded by the number of alternations in a string, and Theorem \ref{Tsofic} which shows that while strings of VC dimension at most $d$ form a shift space (when bi-infinite) this shift space is not sofic when $d>1$. 

\begin{proposition}\label{Pstruct}
Consider a finite binary string $S$ of VC dimension $d \in \mathbb{N}$ with $d>0$.  We assume without loss that $S$ ends in a 1. Suppose $VCdim(S^{([0:-1])})<d$.  Then for some integer $k$, $VCdim(S^{([0:-1])}01^k) = d$.
\end{proposition}
\begin{proof}
Consider $r=\mathfrak{r}(S)$. This point is not in the closed set $\mathfrak{V}_{\leq d-1}$, and therefore there is some open interval $U$ containing $r$ such that every $t\in U$ has $VCdim(t)\geq d$.  For all sufficiently large integers $k$, the real number corresponding to the binary string $S^{([0:-1])}01^k$ will be inside $U$. Therefore, for all sufficiently large $k$, $VCdim(S^{([0:-1])}01^k) \geq d$.  By Lemma \ref{Lpis}, there is some $k^+$ such that $VCdim(S^{([0:-1])}01^{k^+}) = d$.
\end{proof}

\textbf{Example:} The string 011 has VC dimension 2 and 010 has VC dimension 1.  The string 01011 has VC dimension 2 again.  

\begin{proposition} Prepending or postpending a symbol to a binary string increases the VC dimension by at most 1. 
\end{proposition}
\begin{proof}
Like the proof of Lemma \ref{Lpis}.
\end{proof}

\begin{lemma}\label{Lprepend}
Suppose $S$ is a binary string. Then for all $l \in \mathbb{N}$ it is the case that $VCdim(0^lS)\leq VCdim(S)+2$.
\end{lemma}
\begin{proof}
Without loss we may assume that $VCdim(S)$ is finite, say $d\in \mathbb{N}$. Let $S'$ be the reverse of $S$. The mask dimension of $S'$ is the same as the mask dimension of $S$, which is at most $d$.  Let $S'' = S'0^l$, and $S'''$ the reverse of $S''$. Then $S'''=0^lS$, and 

$$VCdim(S)+2 \geq Mdim(S)+2=Mdim(S')+2 \geq VCdim(S')+1 = VCdim(S'')+1$$

$$VCdim(S'') +1 \geq Mdim(S'') + 1= Mdim(S''') + 1 \geq VCdim(S''')=VCdim(0^lS).$$

\end{proof}

The proof technique used in Lemma \ref{Lprepend} can be generalized to prepending or postpending $0^l$ or $1^l$ to any binary string $S$. It may be necessary to take the bitwise complement of a string as part of this process.  Bitwise complements are easily seen to preserve mask dimension. 

\begin{corollary}\label{Cbuild}
Suppose $S$ is a binary string. Then for all $l \in \mathbb{N}$ and $c\in\{0,1\}$ it is the case that $VCdim(c^lS)\leq VCdim(S)+2$ and $VCdim(Sc^l)\leq VCdim(S)+2$.
\end{corollary}
\begin{proof}
The proof of Lemma \ref{Lprepend}, \textit{mutatis mutandis}. 
\end{proof}

Note that $VCdim(0)=0$, and $VCdim(011)=2$, showing that the bound is tight. 

For a binary string $S$, let $\textbf{alt}(S)$ denote the number of maximal contiguous blocks of $1$'s appearing in $S$.  This is a measure of the number of alternations between 0 and 1 that occur in $S$. It may be infinite. 

\begin{theorem}\label{Talt}
If $S$ is a binary string then $VCdim(S) \leq 2\textbf{alt}(S)$.
\end{theorem}
\begin{proof}
If $S$ is a string of 1's then the bound obviously holds, since $VCdim(S)=1$ and $\textbf{alt}(S)=1$. Otherwise, start with the leftmost $0$ occurring in $S$ and build out the rest of the string using Corollary \ref{Cbuild}. Note that adding 0's on the right does not affect VC dimension. 
\end{proof}

When we investigate statements about  $\mathfrak{V}_{\leq d}$ and  $\mathfrak{V}_{< \infty}$ in the remainder of this section, consider that all the strings in both of these sets are without loss members of $2^\mathbb{N}$.  That is, they can be regarded as infinite on the right by extending with zeros if necessary. Adding trailing zeros to a string does not affect VC dimension.  

\begin{corollary}\label{Cshift}
The set of strings $\mathfrak{V}_{<\infty}$ is closed under shifts. 
\end{corollary}
\begin{proof}
This is a consequence of Lemma \ref{Lprepend}.
\end{proof}

This means that $\mathfrak{V}_{<\infty}$ is a subshift.  It is clearly not a subshift of finite type (ie characterized by a finite set of forbidden words).  In fact any finite word occurs in some string of finite VC dimension. Thus $\mathfrak{V}_{<\infty}$ is not closed (or open) in the product topology on $2^{\mathbb{N}}$.  On the other hand  $\mathfrak{V}_{\leq d}$ is not closed under shifts (consider the strings 11 and 011).  However it is topologically closed in $2^{\mathbb{N}}$ with the product topology.  
  Let $\mathfrak{V}_{\leq d}^{\Leftrightarrow}$ denote the set of bi-infinite binary strings of VC dimension at most $d$. 

\begin{proposition}
{$\mathfrak{V}_{\leq d}^{\Leftrightarrow}$ is a shift space.}
\end{proposition}
\begin{proof}
The VC dimension of a string depends only on its substrings.  For bi-infinite strings this set is not affected by shifts, and hence the VC dimension is not affected either.  Thus the set is closed under shifts.  The set $\mathfrak{V}_{\leq d}^{\Leftrightarrow}$ is easily seen to be closed in the product topology on $2^{\mathbb{N}}$.  Therefore it is a shift space.
\end{proof}

Note that both $\mathfrak{V}_{\leq d}$ and $\mathfrak{V}_{\leq d}^{\Leftrightarrow}$ are characterized by omitting finite strings of dimension $d+1$. Thus finite strings of VC dimension $d+1$ form a forbidden list for both sets. We might ask whether $\mathfrak{V}_{\leq d}^{\Leftrightarrow}$ is a shift of finite type (meaning it is characterized by a finite set of forbidden substrings). For $d=0$ this is clearly true, because the forbidden list can simply be $\mathcal{F}=\{1\}$. In Theorem \ref{T2prime} we explicitly show that $\mathfrak{V}_{\leq 1}^{\Leftrightarrow}$ is not a shift of finite type. In fact when $d > 1$, $\mathfrak{V}_{\leq d}^{\Leftrightarrow}$ does not even satisfy a weaker condition, known as being \textit{sofic}. We will not directly define a sofic shift space, but rather use an equivalent condition (see \cite{louidor2013independence}). 

Let $B(\mathfrak{V}_{\leq d}^{\Leftrightarrow})$ denote the set of finite substrings that occur in any element of $\mathfrak{V}_{\leq d}^{\Leftrightarrow}$. 
For a word $w\in B(\mathfrak{V}_{\leq d}^{\Leftrightarrow})$, the \textit{follower set} of $w$ in $\mathfrak{V}_{\leq d}^{\Leftrightarrow}$, denoted by $F_{\mathfrak{V}_{\leq d}^{\Leftrightarrow}}(w)$, is defined
by

$$F_{\mathfrak{V}_{\leq d}^{\Leftrightarrow}}(w) = \{z\in B(\mathfrak{V}_{\leq d}^{\Leftrightarrow}) : wz\in B(\mathfrak{V}_{\leq d}^{\Leftrightarrow})\}.$$

A $\mathbb{Z}$ shift space (such as $\mathfrak{V}_{\leq d}^{\Leftrightarrow})$ is sofic if and only if it has only finitely many follower sets.  We will show that when $d>1$, $\mathfrak{V}_{\leq d}^{\Leftrightarrow}$ has infinitely many distinct follower sets. First we need two technical lemmas. 

The following lemma basically establishes that it is possible for a symbol in a binary string to be so distant that it cannot affect VC dimension.  We use a left infinite string (denoted by $0^\infty$) in the argument. 

\begin{lemma}\label{Lmadness}
Suppose that a binary string $S$ is of the form $0^{\infty}T1$ where 
\begin{enumerate}
\item $VCdim(S)=d+1$ for some $d\in \mathbb{N}$, $d>1$.
\item $VCdim(0^{\infty}T) = d$.
\end{enumerate}
Then there is some $N \in \mathbb{N}$ such that for all $l>N$, $VCdim(0^{\infty}T0^l1)=d$.
\end{lemma}
\begin{proof}
Let $k$ be the length of $T$, and $N=2k$. We claim this is sufficient.  For suppose $l > N$ and $VCdim(0^{\infty}T0^l1)=d+1$.  Then some $A\subseteq \mathbb{N}$ is shattered with $|A|=d+1$. Without loss, $0\in A$.  There is some $A_0 \subseteq A$ that is traced by $0^{\infty}T0^l1$ but not traced by $0^{\infty}T$.  We now argue that $A_0$ contains an element $c>k$.

Because $0^{\infty}T$ is left infinite and $T$ without loss begins with 1, $A_0$ is not a singleton, as all singletons are already traced. If $|A_0|\geq 2$, $\max{A_0} > k$, otherwise $A_0$ would be traced by $0^{\infty}T$. Therefore $\exists c \in A$, $c>k$.  Because $d>1$, there are $a,b \in A$ with $a<b<c$, and $a=0$.

Because $c>k$ and $len(T)=k$, any substring $w$ of $0^{\infty}T0^l1$ tracing a non-singleton set including $c$ must be a suffix. The length of the suffix tracing $\{0,c\}$ must be $c+1$.  But the length of the suffix tracing $\{0,b,c\}$ must also be $c+1$.  Therefore these are the same suffix $\rightarrow\leftarrow$. 

Therefore $VCdim(0^{\infty}T0^l1)=d$.


\end{proof}

The proof shows that actually $0^{2len(T)+1}$ is sufficient, where $0^{\infty}$ appears in the above lemma.

\begin{lemma}\label{Linsert}
Let $S$ a binary string. For a an integer $p>1$, let $S^{[p]}$ be derived from $S$ by replacing each symbol $c$ in $S$ with $0^{p-1}c$. Then $VCdim(S)=VCdim(S^{[p]})$. 
\end{lemma}
\begin{proof}

 It is clear that $VCdim(S^{[p]})\geq VCdim(S)$. We now show $VCdim(S^{[p]})\leq VCdim(S)$.  Suppose $S^{[p]}$ shatters a set of size $k$. Then there is some $A\subseteq \mathbb{N}$ shattered by $S^{[p]}$ with $0\in A$, $|A|=k$ and where $A$ consists only of multiples of $p$.  This $A$ can be shattered even restricting to substrings of the form $S^{[p]([i:j])}$ where $i$ is a multiple of $p$.  In fact if $i \not \equiv 0 \mod p$, the subset of $\mathbb{N}$ corresponding to $S^{[p]([i:j])}$ contains no multiples of $p$. But then $S$ can shatter a set of size $k$ as well. 
\end{proof}

\begin{theorem}\label{Tsofic}
$\mathfrak{V}_{\leq d}^{\Leftrightarrow}$ is not sofic if $d>1$. 
\end{theorem}
\begin{proof}

Let $d>1$. Suppose a binary string $S$ is of the form $0^{2k+1}T1$ where 
the following hold: $k=len(T)$ for $k\in \mathbb{N}$, $VCdim(S)=d+1$, $VCdim(0^{2k+1}T) = d$.

 Define a sequence of numbers $a_1,a_2,\ldots$ recursively by $a_1 = 1$ and $a_{i+1} = 2ka_i+2$. For each $i =1,2,3,\ldots$, define $w_i = 0^{(2k+1)a_i}T^{[a_i]} =  (0^{2k+1}T)^{[a_i]}$.  By Lemma \ref{Linsert}, $VCdim(w_i)=d$ for all $i \in \mathbb{N}\setminus \{0\}$. Therefore $w_i \in B(\mathfrak{V}_{\leq d}^{\Leftrightarrow})$ for all $i$. 
 
 Suppose $i < j$.  We claim $F_{\mathfrak{V}_{\leq d}^{\Leftrightarrow}}(w_i) \neq F_{\mathfrak{V}_{\leq d}^{\Leftrightarrow}}(w_j)$.  In particular we will argue that $0^{a_j-1}1 \in F_{\mathfrak{V}_{\leq d}^{\Leftrightarrow}}(w_i) \setminus F_{\mathfrak{V}_{\leq d}^{\Leftrightarrow}}(w_j)$.  By the definition of $S$ and Lemma \ref{Linsert}, $VCdim(w_j0^{a_j-1}1)=VCdim((0^{2k+1}T1)^{[a_j]})=d+1$.  Therefore $0^{a_j-1}1 \notin F_{\mathfrak{V}_{\leq d}^{\Leftrightarrow}}(w_j)$. 
 
 On the other hand  $w_i = 0^{(2k+1)a_i}T^{[a_i]}$ and $len(T^{[a_i]})=ka_i$.  But $a_j-1 > 2ka_i$.  Therefore, by Lemma \ref{Lmadness}, $$VCdim(w_i0^{a_j-1}1)=VCdim(0^{(2k+1)a_i}T^{[a_i]}0^{a_j-1}1)=d,$$
and $0^{a_j-1}1 \in F_{\mathfrak{V}_{\leq d}^{\Leftrightarrow}}(w_i)$. Thus $F_{\mathfrak{V}_{\leq d}^{\Leftrightarrow}}(w_i) \neq F_{\mathfrak{V}_{\leq d}^{\Leftrightarrow}}(w_j)$ and consequently $\mathfrak{V}_{\leq d}^{\Leftrightarrow}$ is not sofic.





\end{proof}
\section{Prime strings}

Let $d$ be a nonnegative integer and $S$ a binary string. We say that $S$ is $d$-prime for $d \in \mathbb{N}$ if $VCdim(S)=d$ and whenever $S'$ is a proper substring of $S$, $VCdim(S') < d$. In this section we seek to determine the prime strings for various $d$ and analyze the properties of prime strings in general.

A consequence of Lemma \ref{Lmonotone} is that a finite string $S$ is $d$-prime iff the VC dimension of $S$ is $d$ and the dimension decreases when either its leftmost or rightmost symbol is removed.

We now establish the existence of prime strings, and show that all strings contain at least one prime substring of the same VC dimension.

\begin{proposition}\label{Pss}  Let $S$ be a string of VC dimension $d \in \mathbb{N}$.  Then $S$ has a $d$-prime string as a substring.
\end{proposition}

\begin{proof}  We can assume without loss that $S$ is finite.  The following algorithm will produce a prime substring of $S$.  First remove as many symbols from the right of $S$ as possible while preserving the VC dimension.  Then remove as many symbols from the left of the remainder as possible while preserving the VC dimension.  Let $S'$ denote the result of this process. Then $VCdim(S')=d$.  Suppose, by way of contradiction, that $S'$ has a proper substring of VC dimension $d$. Then, by Lemma \ref{Lmonotone}, $S'$ can have a digit removed either from the left or right while preserving the VC dimension. This contradicts the process that produced $S'$ in the first place.
\end{proof}

\begin{lemma}\label{Lpis}:  Let $S$ be a possibly infinite binary string of VC dimension $d+1$ for $d \in \mathbb{N}$.  Then $S$ contains a proper initial substring of VC dimension $d$.
\end{lemma}
\begin{proof}  By Lemma \ref{Lfinite} we can assume without loss that $S$ is finite.  Let $S^+$ be the shortest proper initial substring of $S$ of VC dimension $d+1$. Let $S^-$ be $S^+$ with the rightmost symbol removed.  Then $VCdim(S^-) \leq d$.  We claim that also $VCdim(S^-) \geq d$.

Let $A \subseteq \mathbb{N}$ be a set of size $d+1$ that is shattered by $S^+$.  Let $a = \max(A)$ and $A^- = A\setminus \{a\}$.

We claim that $A^-$ is shattered by $S^-$.  Let $A_0 \subseteq A^-$.  There is some substring $s$ of $S^+$ which cuts out $A_0 \cup \{a\}$ from $A$.

Then $s_a = 1$ and without loss $a$ is the last index in $s$. Let $s^-$ be $s$ without the rightmost symbol, and note that $s^-$ is a substring of $S^-$.  Thus $A_0$ is cut out from $A$ by a substring of $S^-$.  Since $A_0$ was arbitrary, $A$ is shattered by substrings of $S^-$.  Therefore $VCdim(S^-) \geq d$, whence $VCdim(S^-)=d$.
\end{proof}

\begin{proposition}  Let $S$ be a possibly infinite binary string of VC dimension $d+1$ for $d \in \mathbb{N}$. Then $S$ contains a $d$-prime substring.  In particular, every $d+1$-prime string contains a $d$-prime string as a substring.
\end{proposition}
\begin{proof}  First apply Lemma \ref{Lpis} to get a proper initial substring of VC dimension $d$.  Then apply Proposition \ref{Pss} to get a $d$-prime substring.
\end{proof}

Notice that the $d$-prime strings are not closed under complement. For example $ 011 $ is $2$-prime, but $ 100$ is not even VC dimension 2. This is essentially because of the downward closure of the substring concept and the way we associate substrings with subsets in the definition of string VC dimension. This is not the case for mask dimension. Mask dimension is easily seen to be preserved under complementation for any string and reversals for finite strings.

%
%

We now consider prime strings of VC dimension 1,2, and 3.

\subsubsection{VC dimension 1}

There is a unique prime string $S =  1 $.

%
%
%
%

\subsubsection{VC dimension 2}

\begin{theorem}\label{T2prime} Let $k,d$ be nonnegative integers. All 2-prime strings are either of the form
\begin{enumerate}
\item  $10^k10^d1$ where $k < d$, or
\item  $0^{d+1}10^d1$.
\end{enumerate}
\end{theorem}
\begin{proof}
Let $S$ be a 2-prime string.  We may assume that $S$ is finite and ends with 1. By definition, we must have $VCdim(S^{([0:-1])})=1$.  Therefore $S^{([0:-1])}$ is of one of the three forms presented in Lemma \ref{Lvc1}.  Viewed another way, $S$ must result from appending $0^d1$ to one of the forms in Lemma \ref{Lvc1}, for some nonnegative integer $d$.  The forms from Lemma \ref{Lvc1} are as follows. 
\begin{enumerate}
\item $0^a1$, $a\geq 1$
\item $1^a$, $a \geq 1$ 
\item $0^a(10^b)^c1$, $0 \leq a \leq b$, $b\geq 1$, $c\geq 1$ 
\end{enumerate}

We now try appending $0^d1$ to each of these forms in a series of cases (respectively). 

Case 1: Suppose $S=0^a10^d1$ for $a \geq 1$, $d \geq 0$. To achieve $VCdim(S)=2$ it is necessary and sufficient that $d<a$.  Then in order for $S$ to be prime we must have $a=d+1$ (else we could remove a zero on the left while preserving dimension). This gives a 2-prime string of form (2) from the statement of the theorem. 

Case 2: Suppose $S=1^a0^d1$ for $a \geq 1$. To achieve $VCdim(S)=2$, it is necessary and sufficient that $a>1$ and $d>0$. Then for $S$ to be prime, we must have $S=110^d1$, which is form (1) from the statement of the theorem. 

Case 3: Suppose $S=0^a(10^b)^c10^d1$ where $0 \leq a \leq b$, $b\geq 1$, and $c\geq 1$. To achieve $VCdim(S)=2$, it is necessary and sufficient that $d\neq b$. If $d< b$ then it must be that $S=0^{d+1}10^d1$, which is form (2) in the statement of the theorem.  If $d > b$ then we must have $S=10^b10^d1$, which is form (1) in the statement of the theorem.
\end{proof}

It is natural to wonder where $d$-prime strings lie in the hierarchy of formal languages. One could pose similar questions about finite strings of bounded VC dimension in general.  

\begin{corollary}
The language of 2-prime strings is not regular.
\end{corollary}
\begin{proof}
Suppose by way of contradiction that the language is regular.  Then by the Pumping Lemma, all sufficiently long 2-prime strings can be written as a product of words $xyz$ such that $xy^nz$ is still 2-prime for all $n\in \mathbb{N}\setminus\{0\}$ (with $y$ nonempty).  But consider a 2-prime string of the form $S = 0^{d+1}10^d1$, where $d$ is sufficiently large for the Pumping Lemma to apply.  Suppose $S=xyz$.  Then $y$ either does or does not contain a 1. If $y$ does contain a 1, then $xy^2z$ is not 2-prime because it begins with 0 but contains more than two 1's. If $y$ does not contain a 1 then $xy^2z$ is not 2-prime because the 0's are imbalanced. Therefore the language cannot be regular. 
\end{proof}

\subsubsection{VC dimension 3}

The language of 3-prime strings is much more complex than in the case of dimension 2. This is owing to the comparative lack of structure in strings of VC dimension 2, which form an (co)uncountable set (by Theorem \ref{T:uncountable}). In particular, the proof strategy in Theorem \ref{T2prime} cannot be adapted to the case of 3-prime strings, because dimension 2 strings have no simple characterization. 

 The prime strings of VC dimension 2 divide clearly into two infinite families, and have a bounded number of alternations between 0 and 1. However in the case of dimension 3 this is not the case.  In this section we will analyze the shortest 3-prime string $S=0010111$ (found through exhaustive search). This string has an associated family of 3-prime strings, of which it is the simplest member.  These are of the form $00(10)^k 111$, $k \geq 1$.  This family provides examples of 3-prime strings with unboundedly many alternations.  Other examples of 3-prime families  discovered experimentally are presented (without proof) in Table \ref{Ttprime}.
 
\begin{table}
\begin{center}

\begin{tabular}{|l|}
\hline
$S=01111^k0^l10^m1$, where $l \leq k+1, m,k\geq 1$ \\
\hline
$S=0^{2k+1}10^{k-1}10^k10^{k-1}1$, where $k\geq 1$ \\
\hline
$S=1101010101(01)^k0000010000011$, where  $k\geq 0$\\
\hline
\end{tabular}
\caption{Some examples of 3-prime families. }
\label{Ttprime}
\end{center}

\end{table}

In order to simplify the proof of the 3-primeness of $S$ and its family members, we introduce some new terminology.

\begin{definition}
  Given a binary string $S$ of length $n$, the \textit{right rays} of $S$ are the substrings of the form $S^{([i:n])}$ for $0 \leq i < n$.  Let $\mathfrak{R}(S)$ denote the right rays of $s$.
\end{definition}

\begin{definition}  Consider $A\subseteq \mathbb{N}$, consisting of elements $\{a_1 < a_2 < \cdots < a_n\}$.  The \textit{telescope} of $A$ is the family of subsets $\{\{a_1,...,a_i\} : i \leq n\}$.
\end{definition}

\begin{definition}
  Given a substring $s$ of $S$, the telescope of $s$, denoted $\mathcal{T}_s$ is the telescope of the associated set family $n(s) = \{i \in \mathbb{N}: s_i = 1\}$.
\end{definition}

\begin{proposition}
 Given a string S, the set family associated with S as in the definition of VC dimension on strings is

$$\mathfrak{S}=\bigcup_{s \in \mathfrak{R}(S)} \mathcal{T}_s$$

\end{proposition}
\begin{proof}
This is a clear consequence of the definition of the VC dimension of a binary string.
\end{proof}

\begin{proposition}  The family of strings $F_k=00(10)^k 111$, $k \geq 1$ consists only of 3-prime strings.
\end{proposition}
\begin{proof}
There are three things to show.

\begin{enumerate}
\item When $S=00(10)^k 111$, $VCdim(S)=3$ for all $k \geq 1$.
\item When $S=0(10)^k 111$ $VCdim(S)\leq 2$ for all $k \geq 1$.
\item When $S=00(10)^k 11$ $VCdim(S)\leq 2$ for all $k \geq 1$.
\end{enumerate}

In the case of (1) we only need to show that the VC dimension is at least 3, because (3) together with Lemma \ref{Lpis} implies that the VC dimension cannot be 4 or greater.

Proof of (1):

Let $S=00(10)^k 111$ for $k\geq 1$. Suppose $m=len(S)=5+2k$. The subsets of $\mathbb{N}$ corresponding to $\mathfrak{R}(S)$ are described in the following table.  The substrings corresponding to the sets are shown in the column headers.  Each set depends on the row index $i$. 

\begin{table}[h]
\begin{tabular}{|c|c|c|}
\hline 
 $i$& $S^{([-3-2i+1:m])}$ & $S^{([-3-2i:m])}$ \\ 
\hline \hline
0 & $\{0,1\}$ & $\{0,1,2\}$ \\ 
\hline 
1 & $\{1,2,3\}$ & $\{0,2,3,4\}$ \\ 
\hline 
2 & $\{1,3,4,5\}$ & $\{0,2,4,5,6\}$ \\ 
\hline 
3 & $\{1,3,5,6,7\}$ & $\{0,2,4,6,7,8\}$ \\ 
\hline 
4 & $\{1,3,5,7,8,9\}$ & $\{0,2,4,6,8,9,10\}$ \\ 
\hline 
5 & $\{1,3,5,7,9,10,11\}$ & $\{0,2,4,6,8,10,11,12\}$ \\ 
\hline 
$\vdots$ & $\vdots$ & $\vdots$ \\ 
\hline 
$k-1$ & $\{\text{odds to }2k-3\text{ inclusive },2k-2,2k-1\}$ & $\{\text{evens to }2k-2\text{ inclusive },2k-1,2k\}$ \\ 
\hline 
$k$ & $\{\text{odds to }2k-1\text{ inclusive },2k,2k+1\}$  & $\{\text{evens to }2k\text{ inclusive },2k+1,2k+2\}$  \\ 
\hline 
\end{tabular} 
\end{table}

Additionally the following two sets are included, as well as the sets $n(S^{([m:m])})=\emptyset$, and $n(S^{([-1:m])})=\{0\}$.

$$P = n(S^{([1:m])})=\{1,3,5,....,2k+1,2k+2,2k+3\}$$
$$U = n(S^{([0:m])})=\{2,4,6,...,2k+2,2k+3,2k+4\}$$

We will argue that the above set system shatters $A=\{0,2k-1,2k\}$.  In fact $\{0,2l-1,2l\}$ is shattered for $l=1,...,k$ though we omit the proof.

For each subset of $A$ we give the relevant set named above which traces it.  The rows and columns refer to the table above.

\vspace{0.2cm}
\begin{tabular}{|l|l|}
\hline
Subset & Realization \\
\hline \hline 
$\{\}$        &          $S^{([m:m])}$ \\
\hline
$\{0\}$        &         $S^{([-1:m])}$ \\
\hline
$\{2k-1\}$      &        $P$ \\
\hline
$\{2k\}$ &             $U$ \\
\hline               
$\{0,2k-1\}$      &      truncating row $k-1$ column 1 so that $2k-1$ is the maximum element \\
\hline
$\{0,2k\}$      &        row $k$ column 1  \\
\hline
$\{2k-1,2k\}$    &       row $k$ column 0 \\
\hline
$\{0,2k-1,2k\}$   &      row $k-1$ column 1 \\
\hline
\end{tabular}
\vspace{0.2cm}

This concludes the proof of (1).

For (2) and (3) we give only a sketch.  First enumerate the sets realized from right rays.  Then consider a shattered set $A$ of size 3. Work through the various cases concerning whether $A$ contains only odds, only evens, or a mix of odds and evens.  In each case it will be clear that there is some subset of $A$ which is not traced.

\end{proof}

We are uncertain how complex the language of 3-prime strings may be.

\textbf{Question}:  Which level of the hierarchy of formal languages do $d$-prime strings occupy?

\section{Model theory and generalizations}\label{Smt}

The study of binary strings is equivalent to the study of subsets of $\mathbb{N}$ or $\mathbb{Z}$ (if the strings are infinite in both directions.)  The complexity of subsets of integers is a topic that has been deeply investigated in model theory.  These results allow for the quick determination of the finitude of VC dimension for a broad collection of binary strings. In this section we will do two things.  First we will put the results of the previous sections in a model theoretic context and show that many of the results in the previous sections have analogues in structures other than $\mathbb{Z}$ and $\mathbb{N}$. Then we will review results from the model theory literature that establish finitude of VC dimension in a wide assortment of strings. See \cite{marker2006model} for a survey of model theory. 


In the context of model theory a \textit{language} is a set of symbols which denote abstract functions, predicates, and constants.  For example the language of groups is $L=\{+,0\}$ and the language of ordered rings is $L=\{+,\cdot,<\}$.  The equality relation is implicitly included.  If $L$ is a language, then there is an associated set of well formed first order formulas in $L$. A universe of objects together with an interpretation of the elements of $L$ is known as a $L$-structure (or model). For example $N=(\mathbb{N},+,\cdot)$ is a model in the language of rings.  A $L$-model provides a semantics for the $L$-formulas. 

In the presence of a model, formal expressions in a language take on a truth value. If a formal $L$-sentence $\varphi$ is true in an $L$-model $M$, we write $M \models \varphi$, which we read as ``$M$ models $\varphi$". For example $\mathbb{N} \models \forall x \forall y (xy = yx)$ with respect to standard multiplication.  We refer to $M$ both as a model and the universe of the model $M$ by abuse of notation.  The first-order theory of the $L$ structure $M$ is the set of first order $L$-formulas modeled by $M$.

  We say that a language $L'$ is an expansion of a language $L$ if $L \subseteq L'$.   If $\varphi(x_0,x_1,\ldots,x_m)$ is an $L$ formula we can partition its variables such that some are construed as parameters.  When we write $\varphi(\bar{x};\bar{y})$ we mean an $L$ formula with variables $\bar{x}$ of some arity $|\bar{x}|$ and parameter variables $\bar{y}$ with some arity $|\bar{y}|$. If $M$ is an $L$ structure we use $M^{\bar{x}}$ and $M^{\bar{y}}$ to refer to the $|\bar{x}|$-tuples and $|\bar{y}|$-tuples in $M$, respectively. We define, for any $\bar{b} \in M^{\bar{y}}$, $\varphi(M,\bar{b}) = \{\bar{a}\in M^{\bar{x}}: M \models \varphi(\bar{a},\bar{b})\}$.  For any $X \subseteq M^{\bar{x}}$ and $Y \subseteq M^{\bar{y}}$, we have the set system $\mathcal{C}^{\varphi(\bar{x},\bar{y})}_{X,Y} = \{\varphi(X,\bar{b}) : \bar{b} \in Y\}$.  The VC dimension of $\varphi(\bar{x},\bar{y})$ with respect to the theory of $M$ is defined to be the VC dimension of $\mathcal{C}^{\varphi(\bar{x},\bar{y})}_{X,Y}$ where $X=M^{\bar{x}}$ and $Y=M^{\bar{y}}$.  This is first order definable and a property of the theory of $M$.   A first order $L$ theory $T$ is said to be NIP (for ``not the independence property") if every partitioned $L$ formula has finite VC dimension.  A guide to NIP theories can be found in \cite{simon2015guide}. 
  
We now show that what we have called mask dimension on binary strings is simply the VC dimension of certain formulas in Presberger arithmetic, expanded by a predicate. 

 Let $M$ be an $L$-structure where $L$ is an expansion of the language with a single binary operator, denoted $+$.  Let $P(x)$ be a predicate on $M$.  To connect with binary strings, use the convention that $P(a)$ is identified with $1$ if $M \models P(a)$ for $a \in M$, and $P(a)$ is identified with 0 if $M \models \neg P(a)$.  A $d$-mask on $(M,+,P)$ is a set of sequences of the form $\langle P(a_1+t),P(a_2+t),\ldots,P(a_d+t)\rangle$ for fixed $a_1,...,a_d \in M$ as $t$ varies in $M$. A $d$-mask on $(M,+,P)$ is said to be full if $|\{\langle P(a_1+t),P(a_2+t),\ldots,P(a_d+t)\rangle : t \in M\}| = 2^d$. The mask dimension of the predicate $P$ is the maximum $d$ such that a full $d$-mask on $(M,+,P)$ exists.  If there is no maximum, we say that the mask dimension is $\infty$.  
 
If $M=\mathbb{N}$ and $+$ is standard addition, then $P \subseteq \mathbb{N}$ and $P$ can be identified with a binary string.  The mask dimension of this string as defined in Section \ref{Sintro} is what we have defined as the mask dimension of $(M,+,P)$ in the previous paragraph. Note that this is also the VC dimension of $\varphi(x;y) = P(x+y)$ in the theory of $M$.

We can see from this description that all set systems arising from binary strings are self dual, in the sense that $P(x+y)=P(y+x)$.  This is a fairly strong restriction on the kinds of set systems that can arise from binary strings.  

Observe that the entire discussion about mask dimension of a definable set can be applied to more general situations through the notion of the VC dimension of $P(x+y)$.  Below we show that the Cantor set has infinite mask dimension. 

\begin{proposition}  Let $M = (\mathbb{R},+,P)$ and $P(\mathbb{R})$ be the Cantor set.  Then $(\mathbb{R},+,P)$ has infinite mask dimension.
\end{proposition}
\begin{proof}
  An element $r \in \mathbb{R}$ is in the Cantor set iff there is $A \subseteq \mathbb{N}$ such that $r = 2\sum_{a \in A} 3^{-a}$.
  Let $d \in \mathbb{N}$ be given.  
  Consider a $d$-mask of the form $\langle P(a_1+t),P(a_2+t),\ldots,P(a_d+t)\rangle$ where $a_i = 2\cdot 3^{-i}$ for $i =1,2,\ldots,d$.  
  We claim that this $d$-mask is full.  Let $A \subseteq [d]$ be given. 
 Define $t= 2\sum_{a \in A}3^{-a}$.
Now consider $a_i + t = 2(3^{-i}+ \sum_{a \in A} 3^{-a})$.  We have $M \models P(a_i + t) \iff i \notin A$.  This gives a full $d$-mask.  Since $d$ was arbitrary, the Cantor set has infinite mask dimension.
\end{proof}

It was shown in \cite{hieronymi2018interpreting} using model theoretic techniques that an expansion of $(\mathbb{R},+,<)$ that defines a Cantor set is not NIP. The order relation is not necessary (see above) for the standard ternary Cantor set.  The paper however deals with a more abstract definition of a Cantor set, namely a subset of $\mathbb{R}$ that is nonempty, compact, and has neither isolated nor interior points.

It is certain that something equivalent to the following is known, but we cannot find a reference.

\begin{proposition}\label{Pgroup}  Suppose $R=(R,+)$ is a group.  Suppose that $P$ names a subgroup of $R$.  Then the mask dimension of $(R,+,P)$ is one unless $R=P$ in which case it is zero.
\end{proposition}
\begin{proof}
  Let $G$ be the subgroup named by $P$. Let $r_1,r_2 \in R$, and consider a 2-mask $\langle P(r_1+t),P(r_2+t)\rangle$. For this mask to be full, there must be four values of $t$, each realizing a different group membership condition on $r_1+t$ and $r_2+t$.  In particular we must show that there is some $t_{11}$ such that $r_1+t_{11}\in G$ and $r_1+t_{11}\in G$, as well as some $t_{01}$ such that $r_1+ t_{01} \notin G$ and $r_2+t_{01} \in G$. We show that this is not possible.  Without loss $R\setminus G$ is nonempty.
  
    Suppose that there is some $t_{11} \in R$ such that  $r_1+t_{11} = g_1, r_2+t_{11} = g_2$, with $g_1,g_2 \in G$.

Then $g_1-r_1 = g_2 - r_2$, and $g_1-g_2 = r_1 -r_2$.   Therefore $r_1 - r_2 \in G$.

Now suppose that there is some $t_{01} \in R, h \in R \setminus G$, and $g \in G$ such that $r_1 + t_{01} = h$, and $r_2 + t_{01} = g$.

Then $h-g = r_1 - r_2 \in G.$

Therefore $h-g \in G$ and $h-g+g \in G$.  Then $h \in G$ $\rightarrow \leftarrow$.  This shows that the mask dimension of $P$ is less than two.

Suppose that there is some $h \in R\setminus G$.  Let $g \in G$.  Then $g+0 \in G$, and $g+h \notin G$ (since $g+h$ is in the coset $G+h$, disjoint from $G$).  Therefore a full 1-mask exists.

Finally suppose that $R = G$.  Then there is no full 1-mask.
\end{proof}

Proposition \ref{Pgroup} gives many examples of pairs $(M,+,P)$ for which $P$ is of mask dimension 1.

%


The example $(\mathbb{R},+,\mathbb{Q})$ shows in particular that a dense/codense subset of $\mathbb{R}$ can have small complexity. On the other hand a dense/codense subset of $\mathbb{R}$ can have infinite complexity, as the following shows.

\begin{proposition} $(\mathbb{R},+,P)$ has infinite mask dimension where $P$ is the union of sets of the form $\mathbb{Q}+\sum_{i\in A} \pi^i$ for finite subsets $A \subseteq \mathbb{N}$, and $\pi$ a nonalgebraic constant.
\end{proposition}
\begin{proof}
   We identify $P$ with its interpretation $P(\mathbb{R})$.  Because $P \supseteq \mathbb{Q}$, $P$ is dense in $\mathbb{R}$.  Because $P$ is countable it is also codense. Let $d \in \mathbb{N}$ be given.  Let $a_i = \pi^i$ for $i=1,2,...,d$.  Fix any $A \subseteq [d]$ and let $t = \sum_{i \in A} \pi^i$.  Then $a_i + t \in P \iff i \notin A$. This gives a full $d$-mask.
\end{proof}

Note that the definition of a Sidon set generalizes to abelian groups other than $(\mathbb{Z},+)$.  In these more general settings, predicates $P$ which realize near Sidon sets will still have finite mask dimension, essentially by the same argument from Section \ref{Ssidon}.  This is somewhat implicit in results from \cite{VPD2} and \cite{POINT1}. 

\subsubsection{Model theoretic connections}\label{Smodelcon}

In this section we survey the model theory literature for results that give information about classes of binary strings (and generalizations) with finite VC dimension. Often in model theory a goal is to show that a structure has combinatorial properties, such as NIP.  A stronger condition than NIP is \textit{stability}.  Many authors have showed conditions on $P$ such that $(\mathbb{Z},+,P)$ is a stable or NIP structure.  Either of these conclusions implies that every formula, including $P(x+y)$, has finite VC dimension, and hence that the binary string corresponding to $P$ has finite VC dimension.

	There is deep work by a number of authors that examine the groups $G$, subsets $A \subseteq G$, and the VC dimension of the set family $\{gA: g \in G\}$.  This is of course equivalent to studying models $(G,+,A)$ of finite mask dimension where $(G,+)$ is a group (only in multiplicative notation). Examples of work in this vein includes \cite{conant2018structure, conant2018pseudofinite, conant2020approximate}, and \cite{terry2020quantitative}.

Model theorists have been aware for some time that that when $A \subseteq{N}$ is satisfies certain sparsity conditions, $(\mathbb{Z},+,A)$ has finite mask dimension. Some work on this topic includes \cite{conant2019stability,POINT1,conant2018multiplicative}, and \cite{conant2020weakly}. 

Hawthorne has done work on examining the relation between automatic sequences and finite VC dimension \cite{hawthorne2020automata}. 

In \cite{kaplan2017decidability} it is shown (as a consequence of the Green-Tao theorem) that if $P$ names the primes in $\mathbb{N}$, then for all $d$, if $A=[d]$ then for every $A_0 \subseteq A$ there is an arithmetic progression $\{t_0i+t_1: i \in A\}$ such that $P(t_0i+t_1) \iff i\in A_0$. 
Consequently $|\{\langle P(t_0+t_1),P(2t_0+t_1),\ldots,P(dt_0+t_1) \rangle : t_0,t_1 \in \mathbb{N}\} | = 2^d$. This is a kind of ``affine" mask dimension on the set of primes, which is shown to be infinite. 

\textbf{Question:} What is the mask dimension of the binary string with a 1 at index $i$ precisely if $i$ is prime? 

It follows from the work of \cite{POINT1,VPD2} that the sequence with a 1 precisely in indexes corresponding to Fibonacci numbers has finite VC dimension. These authors also first established basically all of the facts in Section \ref{Spowersof2}.  There has been some interesting work relating NIP theories to dynamics on bi-infinite binary strings through the automorphism group \cite{mofidi2018some}.

\bibliographystyle{amsplain}
\small\bibliography{S0168007220300336.bib} 
%


\end{document}